\documentstyle[12pt]{amsart}
\begin{document}

\title[Interpolation and a.e. convergence of integral operators]
{An interpolation theorem related to the a.e. convergence of integral operators}
\author{Alexander Kiselev}
\address{Mathematical Sciences Research Institute, 1000 Centennial Dr., Berkeley, CA 94720, USA}
\email{kiselev@@msri.org}

\begin{abstract}
We show that for integral operators of general form the norm bounds in 
Lorentz spaces imply certain norm bounds for the maximal function. As a consequence, the a.e. convergence for the integral operators on the Lorentz spaces follows from the appropriate norm estimates. 
\end{abstract}

\maketitle

\begin{center}
\section{ Introduction}
\end{center}

In this paper, we study the question of a.e. convergence of integral operators satisfying certain norm estimates in most general setting. The main tool in this studies is an interpolation theorem yielding the norm bounds on the maximal function of the integral operator given norm estimates on the operator itself. 
  The result we show looks so natural and is so general that we hope it may be useful in many situations. In particular, in the Appendix we sketch the application to the spectral analysis of Schr\"odinger operators. It was hard to believe that the result we prove here is not known, but discussions with many specialists in the field convinced me that at least the reference is not simple to find. 

The main result we prove here is as follows.

Let us consider two measure spaces $(X, \mu)$ and $(Y, \nu)$ with 
$\mu,$ $\nu$ being positive measures. Let $A(k,x)$ be a measurable
function on $M \times N,$ integrable over the sets of finite measure in 
$M$ for a.e.~$k \in N.$ Let an operator $T$ be given by
\begin{equation}
(Tf)(k) = \int\limits_{X} A(k,x) f(x)\,d\mu (x). 
\end{equation}
The operator $T$ may be defined originally for all simple functions
or bounded functions of finite support. Consider an extending family
of measurable regions $\Omega_{t} \subset X,$ depending on the real
parameter $t,$ so that $\Omega_{t_{1}} \subset \Omega_{t_{2}}$ if 
$t_{2} > t_{1}.$ Define a maximal function $M_{T},$ corresponding to 
the family $\{ \Omega_{t} \}_{t \in R}$ and the operator $T$ by
\begin{equation}
 (M_{T}f)(k) = \sup_{t} \left|\,\,\int\limits_{\Omega_{t}} A(k,x)f(x)\,dx \right|. 
\end{equation}
We use the notation $\|\,\,\|_{pq}$ for the norm in the Lorentz space $L_{pq}$ or for a conventional quasinorm when the norm does not exist.

Then the following theorem holds: \\

\noindent \bf Theorem. \it Suppose that an operator $T$ given by \rm (1) \it
satisfies the  bounds 
\[ \|Tf\|_{s_{i}\infty} \leq C_{i} \|f\|_{p_{i}1}, \,\,\, i=1,2,\]
where all $s_{i},$ $p_{i}$ are more than or equal to $1,$
$s_{1} > s_{2},$ and $p_{1} \ne p_{2}.$ Then the maximal function 
$M_{T},$ given by \rm (2), \it satisfies 
\begin{equation}
\|(M_{T}f) \|_{s_{r}q} \leq C_{s}(q) \|f\|_{p_{r}q} 
\end{equation}
where $s_{r}^{-1}= (1-r)s_{1}^{-1}+rs_{2}^{-1},$ 
$p_{r}^{-1}= (1-r)p_{1}^{-1}+rp_{2}^{-1},$ $0<r<1,$ and $q$ is any number satisfying $1 \leq q \leq \infty.$ \\

\noindent Remark. \rm The theorem above would turn into just generalized
Marcinkiewicz interpolation theorem if we would change $M_{T}$ to 
$T$ in (3). \\

The motivation for studying this problem comes, in particular, from the applications of a.e. convergence results for certain integral operators to the study of the spectrum of Schr\"odinger operators (continuous as well as discrete) \cite{Kis1}, \cite{Kis2}. We briefly sketch  this relationship in the Appendix, referring to the above works for details. 

We note that the above theorem generalizes some old and well-known results on the a.e. convergence of the orhtogonal series and Fourier integrals, in particular, the results of Paley \cite{Pal} and Zygmund \cite{Zyg}. Of course the theorem does not touch on the question of a.e. convergence at the corner points ($r=0,1$). These questions are very subtle even in the orthogonal series setting. In the case of trigonometric series the question about a.e. convergence of the series with coefficients from $l_{2}(Z)$ is equivalent to the celebrated Luzin's conjecture. Although in this case the a.e. convergence holds by the result of Carleson \cite{Car}, it does not hold in general for orthonormal systems of functions  (see, e.g. \cite{Gar}). In particular, in the case of the multiple Fourier series, C.Fefferman showed that one can choose the extending sets $\Omega_{t}$ so that there is no a.e. convergence on $l_{2}(Z^{n}).$
\cite{Fef}. 

\begin{center}
\section{ Proof of the main result}
\end{center}

\noindent \bf Proof. \rm Fix an arbitrary real number $\tilde{r},$ $0<\tilde{r}<1,$ and take $s_{3},$
$s_{4}$ so that $s_{1}>s_{3}>s_{\tilde{r}}>s_{4}>s_{2}.$ Let 
$p_{3},$ $p_{4}$ correspond to $s_{3},$ $s_{4}$ in a usual way.
Also let us fix some number $q'$ such that  $q'>\max \{p_{3}, p_{4} \}.$

 First, we prove the theorem under an additional assumption that $\mu (\Omega_{t})$ depends continuously on $t.$
 The first step is to decompose the support of $f$ into dyadic pieces and estimate certain auxilliary maximal functions. Such approach was used already by Paley \cite{Pal} in his work on a.e. convergence of the series of orthonormal functions. Let $f$ be a measurable bounded function of finite
support and choose $n$ so that $2^{n-1} \leq \mu({\rm supp}(f)) \leq 2^{n}.$ 
 Let the measurable set $E$ be the
support of the function $f:$ $E= \{ x| \;\, |f(x)|>0 \}.$ For every
integer $m<n,$ we consider a partition of the set $E$ into the dyadic pieces of size $2^{m}$ in a following
way:
\[ E_{m,l} = (\Omega_{t_{m,l+1}} \setminus \Omega_{t_{m,l}}) \cap E, \]
where $t_{m,l}$ is defined by a condition
\[ t_{m,l} = \inf \{ t | \mu(\Omega_{t} \cap E)  = 2^{m}l.\} \]
The value $t_{m,l}$ is well-defined at least for $l \leq 2^{n-m-1}$
because $\mu (\Omega_{t})$ depends on $t$ continuously.
The number of the sets $E_{m,l}$ is between $2^{n-m}$ and $2^{n-m-1}.$
For notational convenience, we will assume that this number is always
$2^{n-m}$ and will define the missing $E_{m,l}$ as empty sets.

Let us define the functions $M_{m,l}f$ and $M_{m}f$ by
\[ M_{m,l}f(k) =  \,\int\limits_{E_{m,l}} A(k,x)f(x)\,dx . \]
and
\[ M_{m}f(k) = \sup_{l} |M_{m,l}f(k) | \]
Considering a dyadic development of every real number $N,$ it is easy to see that
\begin{equation}
(M_{T}f)(k) \leq \sum\limits^{n}_{m=-\infty} M_{m}f(k).
\end{equation}
Indeed, suppose that for a given value of $k,$ the supremum in (2) is reached
when the value of the parameter $t$ is equal to $N$ (clearly if $f$ has finite support, the
supremum is reached for some value of $t$). Define the real number $a$ by
$a= \mu(E \cap \Omega_{N})$ and consider the dyadic development of $a:$
\[ a= \sum_{j=-\infty}^{n} a_{j}2^{j},\] where $a_{j}$ is equal
to $0$ or $1$ for every $j.$ Then by construction, we can find disjoint
sets $E_{m,l},$ at most one for each value of $m,$ so that 
\[ \mu\left((E\cap \Omega_{N}) \setminus (\cup_{m} E_{m,l})\right)=0. \]
 In fact, for each $m$ the set $E_{m,l}$ belongs to
the union iff $a_{m}=1.$ The corresponding value of $l$ then may be found
by the formula $l= \sum^{n}_{j=m+1} a_{j} 2^{j-m}.$
 By the generalized Marcinkiewicz interpolation theorem (see, e.g., \cite{BL} or \cite{St}) we have
\[ \|Tf\|_{s_{r}q'} \leq C_{r}(q') \|f\|_{p_{r}q'}, \,\,0<r<1. \]
Recall that the quasinorm $\|\cdot\|^{*}$ in Lorentz spaces may be defined by
\[ \|f\|^{*}_{pq} = \left( \frac{q}{p} \int\limits^{\infty}_{0} |f^{*}(t)|^{q}
t^{\frac{q}{p}-1}\,dt\right) ^{\frac{1}{q}}. \]
Let us denote by $\chi (E)$ the characteristic function of the measurable set $E.$
Using the equivalence of the quasinorm $\|\cdot \|^{*}_{pq}$ and norm $\|\cdot\|_{pq}$ for $p>1$  (note that $p_{i}$ are automatically more than $1$), we see that in particular
\[ \|M_{m,l}f(k)\|^{s_{i}q'}_{q'}= \|T(f \chi (E_{m,l}))(k) \|_{p_{i}q'}^{q'} \leq
C_{i}(q')^{q'}\left( \frac{q'}{p_{i}} \int\limits_{0}^{\infty} | (f\chi (E_{m,l}))^{*}
|^{q'}t^{\frac{q'}{p_{i}}-1} \,dt \right)\leq \]
\[ \leq C_{i}(q')^{q'}2^{m(\frac{q'}{p_{i}}-1)}\frac{q'}{p_{i}} \| f\chi (E_{m,l})(t)\|^{q'}_{q'}, \,\, i=1,2.\]
Note that 
\[ \| M_{m}f\|_{s_{i}q'}^{q'} \leq \sum\limits_{l=1}^{2^{n-m}}
\| M_{ml}f \|_{s_{i}q'}^{q'}. \]
Hence summing over $l$ we obtain
\[ \|M_{m}f\|_{s_{i}q'}^{q'} \leq C_{i}(q')^{q'}\frac{q'}{p_{i}} 2^{m(\frac{q'}{p_{i}}-1)}\|f\|_{q'}^{q'}. \]
By (4), we have that
\begin{equation}
\| M_{T}f \|_{s_{i}q'} \leq C_{i}(q') \left( \frac{q'}{p_{i}} \right)^{\frac{1}{q'}} \|
f \|_{q'} \sum\limits_{m=-\infty}^{n} 2^{m (\frac{1}{p_{i}}-\frac{1}{q'})}=
B_{i}(q') 2^{n(\frac{1}{p_{i}}-\frac{1}{q'})}\|f\|_{q'}.
\end{equation}
In the last inequality we denoted by $B_{i}(q')$ a new constant, which only depends on $i$ and $q'.$
Now we note that in a particular case when $f$ is a characteristic function
of a set, $f=\chi (E),$ (5) means
\begin{equation}
 \| M_{T}\chi (E) \|_{s_{i}q'} \leq B_{i}(q') 2^{n(\frac{1}{p_{i}}-\frac{1}{q'})}2^{\frac{n}{q'}}=
B_{i}(q') 2^{\frac{n}{p_{i}}} \leq 2^{\frac{1}{p_{i}}}B_{i}(q') \| \chi (E) \|_{p_{i}1}. 
\end{equation}
It is easy to verify that the operator $M_{T}f,$ defined originally on the measurable
bounded functions of compact support, is a sublinear operator. 
We remind that it means 
\[ |(M_{T}(f_{1}+f_{2}))(k)|\leq |(M_{T}f_{1})(k)|+|(M_{T}f_{2})(k)|, \,\,\,
|(M_{T}(bf))(k)|=|b||(M_{T}f)(k)| \]
for every scalar $b.$
It is  a well-known and simple to check fact that from the inequality (6) for sublinear operator it follows that
$\|Mg\|_{s_{i}q'} \leq C\|g \|_{p_{i}1}$ holds for all simple functions (see \cite{St}) and hence by simple limiting
argument for all bounded functions of finite support. Interpolating, we obtain that $\|Mf\|_{s_{\tilde{r}}q} \leq C_{\tilde{r}}(q)\|f\|_{p_{\tilde{r}}q}$ for every function $f$ bounded
and of compact support.
 It is straightforward to see
that this relation is then extended to all functions $f \in L_{p_{\tilde{r}}q}.$ This completes the proof under the assumption that $\mu (\Omega_{t})$ is continuous.

Consider now the general case. The technical problem we face here
is that now it is not so easy to apply the dyadic decomposition of the ${\rm supp}(f)=E$ into coherent pieces: the jumps of the monotone function $\mu (E \cap \Omega_{t})$ may in general pose obstacles to that. We will handle this problem
by constructing an auxilliary measure space and an auxilliary operator.
We prove the estimates for the maximal function of this auxilliary operator,
and then show that from these estimates follows the result for the original
problem. 

At every value of the parameter $t$ the function $\mu (\Omega_{t}),$ as a
monotone function, has limits from the left and from the right.  Denote by 
$t_{\mp}$ the value of the jump on the left and on the right respectively:
$t_{+}=\mu(\Omega_{t+0} \setminus \Omega_{t})$ and $t_{-}=\mu (\Omega_{t}
\setminus \Omega_{t-0}).$ The set of values of $t$ where any jump may occur
is clearly at most countable. Let the sequence $\{t_{n}\}_{n=1}^{\infty}$ 
denote these points. To each $t_{n}$ corresponds the size of the jump of $\mu
(\Omega_{t})$ at this point, $y_{n}.$ If for some values of $t$ we have both $t_{-}$ and 
$t_{+}$ nonzero, this value of $t$ is encountered in the above sequence 
twice: say, $t_{n}=t_{n+1}=t$ and $y_{n}=t_{-},$ $y_{n+1}=t_{+}.$ Hence, for each member of the sequence $t_{n}$ only one of the two possible jumps takes place. We denote 
by $\Gamma_{t_{n}}$ the portion of the measure space corresponding to $t_{n}:$
$\Gamma_{t_{n}} = \Omega_{t_{n}+0} \setminus \Omega_{t_{n}}$ if we have a 
jump on the right at $t_{n}$ and $\Gamma_{t_{n}}= \Omega_{t_{n}} \setminus
\Omega_{t_{n}-0}$ otherwise. 

Let us consider the following auxilliary measure space $\tilde{X}$ which we biuld out of $X.$ At each 
value of $t_{n}$ we replace $\Gamma_{t_{n}} \subset X$ by 
\[ B_{t_{n}} = \Gamma_{t_{n}} \times [0, 1] \]
with the structure of the product measure space. Hence
\[ \tilde{X} = \left( X \setminus (\cup_{n} \Gamma_{t_{n}})\right) \cup\left(\cup_{m}B_{t_{m}} \right). \]
 The measure $\tilde{\mu}$ on $\tilde{X}$ coincides with measure $\mu$ on
the measurable set 
\[ X_{0} = X \setminus \left( \cup_{n} \Gamma_{t_{n}} \right), \]
while on $B_{t_{n}}$ the measure $\tilde{\mu}$ equals the product measure
$\mu \times dx$ ($dx$ being a Lebesgue measure on $[0,1]$). Furthermore, we
 let the kernel $\tilde{A}(k,x),$ 
defined on $Y \times \tilde{X},$  be equal to $A(k,x)$ for all $k$ when $x \in M_{0}$ and $\tilde{A}(k, x, y) =A(k, x)$ for all $k \in Y$ and all 
$y \in [0,1]$ if $x$ belongs to $\Gamma_{t_{n}}$ for some $n.$ 
Define an integral operator $\tilde{T}$ by
\[ (\tilde{T}f)(k)= \int\limits_{\tilde{X}}\tilde{A}(k,x)f(x)\,dx.\]
Next, define a family $\tilde{\Omega}_{u}$ of the extending measurable sets 
in $\tilde{X}.$ We construct $\tilde{\Omega}_{u}$ so that $\tilde{\mu}
(\tilde{\Omega}_{u})=u.$ Let $t_{0}(u) = \sup_{t}\{ t| \mu(\Omega_{t})
\leq u \}.$ If $t_{0}(u) \ne t_{n}$ for any $n,$ we let
\[ \tilde{\Omega}_{u} = \left( \Omega_{t_{0}(u)} \cap X_{0} \right) \bigcap
\left( \cup_{t_{m}<t_{0}(u)} B_{t_{m}} \right). \]
In this case $\tilde{\mu}(\tilde{\Omega}_{u})=u,$ since $\mu (\Omega_{t})$
is continuous at every $t \notin \{t_{n} \}_{n=1}^{\infty}.$
Suppose now that $t_{0}(u) = t_{n}.$ If at $t_{n}$ we have a jump on the left
($ \mu(\Omega_{t_{n}} \setminus \Omega_{t_{n}-0})>0$) we let 
\[ \tilde{\Omega}_{u} = \left(\Omega_{t_{n}-0}\cap X_{0}\right) \cup \left(\cup_{t_{m}<t_{n}}
\Gamma_{t_{m}} \times [0,1]\right)\cup \left( \Gamma_{t_{n}} \times [0,
\frac{u-\mu (\Omega_{t_{n}-0})}{\mu(\Gamma_{t_{n}})}]\right). \]
Otherwise, if the jump is on the right ($\mu (\Omega_{t_{n}+0} \setminus 
\Omega_{t_{n}})>0$), we define 
\[ \tilde{\Omega}_{u} = \left(\Omega_{t_{n}}\cap X_{0}\right) \cup \left(\cup_{t_{m}<t_{n}}
\Gamma_{t_{m}} \times [0,1]\right)\cup \left( \Gamma_{t_{n}} \times [0,
\frac{u-\mu (\Omega_{t_{n}})}{\mu(\Gamma_{t_{n}})}]\right), \]
essentially just changing $t_{n}-0$ to $t_{n}.$ Defined this way,
$\tilde{\Omega}_{u}$ constitutes a measurable, extending family of sets
such that $\tilde{\mu}(\tilde{\Omega}_{u})=u$ and, in particular, is
continuous. Define the maximal function $\tilde{X},$ corresponding to 
the operator $\tilde{T}$ and the family $\tilde{\Omega}_{u}:$ 
\[ (\tilde{M}_{\tilde{T}}f)(k) = \sup_{u} \left|\,\, \int\limits_{\tilde{\Omega_{u}}} \tilde{A}(k,x) f(x)\,d\mu(x) \right|. \]
We have \\
\bf Lemma. \it For every $r,$ $1>r>0,$ the operator $\tilde{T}$ satisfies
the norm bounds
\[ \|\tilde{T}f \|_{s_{r}p_{r}} \leq C(r) \|f\|_{p_{r}}. \]

\noindent \bf Proof. \rm Indeed,
\[ \tilde{T}f(k) = \int\limits_{X_{0}} A(k, x)f(x)\,d\mu(x) + \sum\limits_{n}
\int\limits_{\Gamma_{t_{n}}}d\mu(x) \int\limits_{0}^{1} \tilde{A} (k,x,y)f(x,y)\,dy = \]
\[ = \int\limits_{X_{0}} A(k,x)f(x)d\mu(x) + \sum\limits_{n} \int\limits_{\Gamma_{t_{n}}} d\mu(x) A(k,x) \int\limits_{0}^{1} f(x,y)\,dy \]
by the definition of the kernel $\tilde{A}(k,x,y).$ Because of the norm 
bounds on the operator $T$ and interpolation between them, we obtain
\[ \|\tilde{T}f\|_{L_{s_{r}p_{r}}(Y)} \leq C(r) \left( \,\int\limits_{X_{0}}|f(x)|^{p_{r}}d\mu(x)+ \sum\limits_{n}\int\limits_{\Gamma_{t_{n}}}d\mu(x) \left| \int\limits_{0}^{1}f(x,y)\,dy \right|^{p_{r}}\right)^{\frac{1}{p_{r}}} \leq \]
\[ \leq C(r) \left( \int\limits_{X_{0}} |f(x)|^{p_{r}} d\mu (x) + \sum\limits_{n} \int\limits_{\Gamma_{t_{n}}} d\mu(x) \int\limits_{0}^{1}|f(x,y)|^{p_{r}}\,dy\right)^{\frac{1}{p_{r}}} \leq C(r) \|f\|_{L_{p_{r}}(\tilde{X})}. \]
We used Jensen inequality in the second step. $\Box$ \\

Now let us choose $r_{1},$ $r_{2}$ so that $s_{r_{1}}=s_{3}$ and $s_{r_{2}}=s_{4}.$
By the lemma, we have the bounds 
\[ \|\tilde{T}f\|_{s_{i}\infty} \leq C(i)\|f\|_{p_{i}}, \,\,i=3,4. \]
Since the function $\tilde{\mu}(\tilde{\Omega}_{u})$ is continuous in $u,$
we can infer from the first part of the proof that the bound
\[ \|\tilde{M}_{\tilde{T}} \|_{s_{\tilde{r}}q} \leq C_{\tilde{r}}(q)\|f\|_{p_{\tilde{r}}q} \]
hold for every $q \in [1, \infty].$ But for every $t,$
we have 
\[ \int\limits_{\Omega_{t}} A(k,x)f(x)d\mu(x) = \int\limits_{\tilde{\Omega}_{u(t)}}\tilde{A}(k,y)\tilde{f}(y)\,d\tilde{\mu}(y). \]
Here $u(t)$ is such that $\mu (\Omega_{t})=u(t)$ and $\tilde{f}$ is a function 
which we define in a clear way.  
Namely, if $y=x \in X_{0}$ we let $\tilde{f}(y)=f(x)$ and if $y=(x_{1},y_{1}) \in \Gamma_{t_{n}} \times [0,1],$ we let $\tilde{f}(x_{1},y_{1})=f(x_{1}).$ It is easy to see that $\|\tilde{f}\|_{L_{pq}(\tilde{X})}=\|f\|_{L_{pq}(X)}$ for all $p,$ $q$
because the distribution functions of $\tilde{f}$ and $f$ coincide.
Therefore, 
\[ (M_{T}f)(k) = \sup_{t} \left| \,\,\int\limits_{\Omega_{t}}A(k,x)f(x)d\mu(x)\right| \leq \sup_{u}\left|\,\, \int\limits_{\tilde{\Omega}_{u}} \tilde{A} (k,y)\tilde{f}(y)d\tilde{\mu}(y) \right| = \tilde{M}_{\tilde{T}}\tilde{f}(k) \]
and hence 
\[ \|M_{T}f\|_{L_{s_{\tilde{r}}q}(X)} \leq \|\tilde{M}_{\tilde{T}}\tilde{f}\|_{L_{s_{\tilde{r}}q}(\tilde{X})} \leq C_{\tilde{r}}(q) \|\tilde{f}\|_{L_{p_{\tilde{r}}q}(\tilde{X})} = C_{\tilde{r}}(q)\|f\|_{L_{p_{\tilde{r}}q}(X)}. \]
Hence we have shown the bound (3) in the general case and since the number $\tilde{r}$ is arbitrary between $0$ and $1,$ the proof of the theorem is now complete. $\Box$ \\

  We may obtain the Zygmund's and Paley's theorems by specifying the measure spaces and the kernels. Let $X$ and $Y$ to be real lines equipped with Lebesgue measure and $A(k,x)=\exp (ikx).$ The $L_{1}-L_{\infty}$ and $L_{2}-L_{2}$ bounds are obvious and we get the Zygmunds theorem, which says that the Fourier integral converges a.e. for 
functions $f$ from $L_{p},$ $1 \leq p <2$ and claims the corresponding estimate
(3) for the maximal function. 

Next let $X=Z^{1},$ $Y=(a,b)$ and $A(k,n)=\phi_{n}(k),$ where $\{ \phi_{n}(k) \}_{n=1}^{\infty}$ is an orthonormal uniformly bounded system in $L_{2}(a,b).$ Again the same bounds as above hold and we get the Paley's theorem which says that the series of orhtonormal functions $\sum_{n}c_{n}\phi_{n}(x)$ converges a.e. if $c_{n} \in l_{p},$ $1 \leq p <2.$

The last example we would like to present here is that of a pseudodifferential operator
in $R^{n}$ with bounded symbol $a(k,x):$
\begin{equation}
 (Tf)(k)= \int\limits_{R^{n}} \exp (ikx) a(k,x) f(x)\,dx .
\end{equation}
If the symbol $a(k,x)$ is from the class for which an $L_{2}-L_{2}$ estimate is valid, we have that the integral defining (7) convergres a.e.~$k$ no matter which system of extending regions $\Omega_{t}$ we take. This is the situation which appears (with $n=1$) in some of the applications to Schr\"odinger operators \cite{Kis1}.

\begin{center}
 \section*{Appendix}
\end{center}

In \cite{Kis1}, a new approach to the investigation of the stability of the absolutely continuous spectrum of one-dimensional Schr\"odinger operators under slowly decreasing pertubations is developed. The absolutely continuous spectrum
 corresponds to the infinite motion of the quantum particle. The question of the preservation of the absolutely continuous spectrum under decaying perturbations means, roughly speaking, finding out which local, i.e. decaying, pertubations are not strong enough to destroy this infinite motion. It turns out that in many cases the problem of spectral analysis may be reduced to a problem of finding asymptotics of the solutions of a certain ODE system, depending on the parameter $\lambda$, for a.e. values of this parameter ($\lambda$ has a meaning of energy for Schr\"odinger equation). The 
latter problem may be tackled with asymptotic integration methods. The first  step is an introduction of a transformation of the system, containing, in a simplest case, a function 
\begin{equation}
 q(x, \lambda) =  \int\limits_{x}^{\infty} \theta^{2}(y, \lambda)V(y)\,dy. 
\end{equation}
Here $\theta (x, \lambda)$ is a solution of the unperturbed Schr\"odinger equation, so that  
\[ \left(-\frac{d^{2}}{dx^{2}}+U\right)\theta =\lambda\theta \]
and $V(x)$ is a  perturbation. It turns out that if the function $q(x, \lambda)$ may be defined and decays sufficiently fast for the values of $\lambda$ from the support of the absolutely continuous spectrum of the unperturbed operator $(-\frac{d^{2}}{dx^{2}}+U)$, then this absolutely continuous spectrum is preserved under perurbation by $V.$ For example, the simplest condition under which the a.c. spectrum is preserved is that $q(x, \lambda)V(x) \in L^{1}$ for a.e. $\lambda$ from the support of the a.c. spectrum of the unperturbed operator. The interesting case is when $V$ is not absolutely integrable; otherwise the stability is known for a long time and may be proven in a very simple way. However, although there has been a considerable attention to the subject, until recently no other general classes preserving the a.c. spectrum,  given in terms of the rate of decay, were known even in the simplest situations (such as $U=0$). 

Thus one is led to the question of studying a.e. convergence and rate of decay estimates for the integral operators like (8). The paper \cite{Kis1} contains a simpler and less general version of the theorem we proved here. That result suffices for the applications of that paper, in particular
producing a new general class of potentials preserving the a.c. spectrum of the free ($U=0$) and periodic ($U$ periodic) Schr\"odinger operators: all potentials $V$ satisfying $|V(x)| \leq C(1+x)^{-\frac{2}{3}-\epsilon}.$ The more general version we present here, however,  may also be  applied to studying discrete
Schr\"odinger operators \cite{Kis2}.

\begin{center}
\section*{Acknowledgement}
\end{center}

I am very grateful to Prof. S.~Semmes and Prof. G.~Pisier for stimulating discussions. I gratefully acknowledge hospitality of IHES, where part of this work was done. Research at MSRI supported in part by NSF grant DMS 9022140.


\begin{thebibliography}{99}

\bibitem{BL} J.~Bergh and J.~L\"ofstr\"om, \it Interploation Spaces: An Introduction, \rm Sprinder-Verlag, Berlin Heidelberg 1976.

\bibitem{Car} L.~Carleson, \it On convergence and growth of partial sums of Fourier series, \rm Acta Math. {\bf 116} (1966), 135--157.

\bibitem{Fef} C.Fefferman, \it On the divergence of multiple Fourier series, \rm Bull. of the Amer. Math. Soc. {\bf 77} (1976), 87--88.

\bibitem{Gar} A.~Garsia, \it Topics in a.e. Convergence, \rm Markham Pub. Company, Chicago 1970.

\bibitem{Kis1} A.~Kiselev, \it Preservation of the absolutely continuous spectrum of Schr\"odinger equation under perturbations by slowly decreasing potentials and a.e. convergence of integral operators, \rm preprint.

\bibitem{Kis2} A.~Kiselev, \it Stability of the absolutely continuous spectrum of Jacobi matrices under slowly decaying perturbations, \rm in preparation.

\bibitem{Pal} R.E.A.C.~Paley, \it Some theorems on orthonormal functions, \rm Studia Math. {\bf 3} (1931), 226--245.

\bibitem{St} E.M.~Stein and G.~Weiss, \it Introduction to Fourier analysis on Euclidean spaces, \rm Princeton Univ. Press, Princeton 1971.

\bibitem{Zyg} A.~Zygmund, \it A remark on Fourier transforms, \rm Proc. Camb. Phil. Soc. {\bf 32} (1936), 321--327.

\end{thebibliography}
\end{document}